\documentclass[12pt,leqno]{article}

\usepackage{amsfonts}
\usepackage{amssymb}

\newcommand{\ra}{\rightarrow}

\newcommand{\cO}{{\cal O}}

\newcommand{\ZZ}{{\mathbb Z}}
\newcommand{\CC}{{\mathbb C}}

\newcommand{\NN}{{\mathbb N}}

\newcommand{\PP}{{\mathbb P}}

\begin{document}
\begin{center} {\bf\Large Faithful Action on the Space of Global\\
\vspace*{0.2cm} Differentials of an Algebraic Curve} \end{center}
\begin{center}\sc Bernhard K\"ock \end{center}

\bigskip

\begin{quote}
\footnotesize {\bf Abstract.} Given a faithful action of a
finite group on an algebraic curve of genus at least 2, we
prove that the induced action on the space of global
holomorphic differentials is faithful as well, except in the
following very special case: the given action is not tame, the
genus of the quotient curve is 0 and the characteristic of the
base field is 2.


\end{quote}

\bigskip

Let $X$ be a connected smooth projective algebraic curve over
an algebraically closed field $k$ equipped with a faithful
action of a finite group $G$ of order $n$. Then $G$ also acts
on the vector space $H^0(X, \Omega_X)$ of global holomorphic
differentials on $X$. A widely studied problem is to determine
the structure of $H^0(X,\Omega_X)$ as module over the group
ring $k[G]$. It goes back to the Chevalley-Weil when $k = \CC$,
see \cite{CW}. If the canonical projection $\pi: X \rightarrow
Y$ from $X$ to the quotient curve $Y = X/G$ is tamely ramified,
a fairly explicit answer to this problem has been given in 1986
by Kani in \cite{Ka}. For more recent answers to (related)
questions in more general situations the reader is referred to
the papers \cite{Bo} and \cite{FWK}. In the case of arbitrary
wild ramification the explicit calculation of the
$k[G]$-isomorphism class of $H^0(X,\Omega_X)$ is still an open
problem.

This note is concerned with the weaker question whether $G$
acts faithfully on the space~$H^0(X,\Omega_X)$. We give the
following answer to this question. Let $g_X$ and $g_Y$ denote
the genus of $X$ and $Y$, respectively, and let $p$ denote the
characteristic of $k$.

{\bf Theorem.} {\em Suppose that $G$ does not act faithfully on
$H^0(X,\Omega_X)$. Then $g_X \in\{0,1\}$ or all of the
following three conditions hold:\\
(i) The projection $\pi$ is not tamely ramified.\\
(ii) $g_Y =0$.\\
(iii) $p=2$. }

The proof of this theorem will be given after the proof of
Proposition~1 below.

For trivial reasons the group $G$ does normally not act
faithfully on $H^0(X,\Omega_X)$ if $g_X \in \{0,1\}$, see
parts~(a) and~(b) of the following example. That we also cannot
expect the action of $G$ on $H^0(X,\Omega_X)$ to be faithful in
the (very special) situation described by conditions~(i), (ii)
and (iii) of the above theorem, is explained in part~(c) of the
following example and in Proposition~2 below.

{\bf Example.}\\
(a) If $g_X = 0$ and $G$ is not the trivial group then $G$ does
obviously not act faithfully on $H^0(X,\Omega_X) = \{0\}$.\\
(b) If $g_X =1 $ (that is if $X$ is an elliptic curve) and if
$G$ is a finite subgroup of $X(k)$ acting on $X$ by translation
then $G$ leaves invariant any global non-vanishing holomorphic
differential and hence $G$ acts trivially on
$H^0(X,\Omega_X)$.\\
(c) Let $p=2$ and let $r$ be an odd natural number. Let
$k(x,y)$ be the cyclic field extension of the rational function
field $k(x)$ of degree~2 given by the Artin-Schreier equation
$y^2 - y = x^r$. Let $\pi: X \rightarrow \PP^1_k$ be the
corresponding cover of nonsingular curves over $k$. Then $\pi$
is not tamely ramified (see Example~2.5 on p.~1095 in
\cite{Ko}) and the Galois group $G = \ZZ/2\ZZ$ acts trivially
on $H^0(X,\Omega_X)$. This follows from Proposition~2 below.

The next lemma is crucial for the proof of Proposition~1 which
in turn is the main idea for the proof of our theorem. We begin
by introducing some notations. For any $G$-invariant divisor
$D$ on $X$ let $\cO_X(D)$ denote the corresponding equivariant
invertible $\cO_X$-module, as usual. Furthermore let
$\pi_*^G(\cO_X(D))$ denote the subsheaf of the direct image
$\pi_*(\cO_X(D))$ fixed by the obvious action of $G$ on
$\pi_*(\cO_X(D))$ and let $\left\lfloor \frac{\pi_*(D)}{n}
\right \rfloor$ denote the divisor on $Y$ obtained from the
push-forward $\pi_*(D)$ by replacing the coefficient $m_Q$ of
$Q$ in $\pi_*(D)$ with the integral part $\left \lfloor
\frac{m_Q}{n} \right \rfloor$ of $\frac{m_Q}{n}$ for every $Q
\in Y$. The function fields of $X$ and $Y$ are denoted by
$K(X)$ and $K(Y)$, respectively. Finally, for any $P \in X$,
let $e_P$ denote the ramification index of $\pi$ at $P$ and let
$\textrm{ord}_P$ and $\textrm{ord}_Q$ denote the respective
valuations of $K(X)$ and $K(Y)$ at $P$ and $Q:=\pi(P)$.

{\bf Lemma.} {\em Let $D$ be a $G$-invariant divisor on $X$.
Then the sheaves $\pi_*^G(\cO_X(D))$ and
$\cO_Y\left(\left\lfloor \frac{\pi_*(D)}{n}\right
\rfloor\right)$ are equal as subsheaves of the constant sheaf
$K(Y)$ on $Y$. In particular the sheaf $\pi_*^G(\cO_X(D))$ is
an invertible $\cO_Y$-module. }

For the reader's convenience we include a proof of this lemma
although it may already exist in the literature.

{\em Proof.} For every open subset $V$ of $Y$ we have
\[\pi_*^G(\cO_X(D))(V) = \cO_X(D) (\pi^{-1}(V))^G \subseteq K(X)^G = K(Y).\]
In particular both sheaves are subsheaves of the constant sheaf
$K(Y)$ as stated. It therefore suffices to check that their
stalks are equal. Let $Q \in Y$, let $P \in \pi^{-1}(Q)$ and
let $n_P$ denote the coefficient of $D$ at $P$. Then we have
\begin{eqnarray*}
\lefteqn{\pi_*^G\left(\cO_X(D)\right)_Q = \cO_X(D)_P \cap
K(Y)}\\
&=& \left\{f \in K(Y): \textrm{ord}_P(f) \ge -n_P\right\}\\
&=& \left\{f \in K(Y): \textrm{ord}_Q(f) \ge - \frac{n_P}{e_P}\right\}\\
&=& \left\{ f \in K(Y): \textrm{ord}_Q(f) \ge - \left\lfloor
\frac{n_P}{e_P} \right\rfloor \right\}\\
&=& \cO_Y\left(\left\lfloor \frac{\pi_*(D)}{n} \right\rfloor
\right)_Q,
\end{eqnarray*}
as desired. \hfill $\Box$

Let $R := \sum_{P\in X} \textrm{dim}_k (\Omega_{X/Y}) [P]$
denote the ramification divisor of $\pi$. The following
proposition computes the dimension of the subspace of
$H^0(X,\Omega_X)$ fixed by $G$.

{\bf Proposition~1.}
\[\dim_k \left(H^0(X,\Omega_X)^G\right) = \left\{
\begin{array}{ll}
g_Y & \textrm{if } \textrm{deg} \left\lfloor \frac{\pi_*(R)}{n} \right\rfloor = 0\\
\\
g_Y-1 + \textrm{deg}\left\lfloor \frac{\pi_*(R)}{n} \right\rfloor &
\textrm{if } \textrm{deg}\left\lfloor \frac{\pi_*(R)}{n} \right\rfloor > 0.
\end{array}\right.\]

{\em Proof.} Let $K_X$ be a $G$-invariant canonical divisor on
$X$, that is we have an equivariant isomorphism $\cO_X(K_X)
\cong \Omega_X$. Let the divisor $K_Y$ on $Y$ be defined by the
equality $\pi^*(\Omega_Y) = \cO_X(\pi^*(K_Y))$ of subsheaves of
$\cO_X(K_X)$. Note that we consider $\pi^*(\Omega_Y)$ as a
subsheaf of $\Omega_X \cong \cO_X(K_X)$ and that we have a
short exact sequence
\[0 \rightarrow \pi^* \Omega_Y \rightarrow \Omega_X \rightarrow
\Omega_{X/Y} \rightarrow 0. \]
 In particular we have $K_X = \pi^* K_Y + R$ and hence
\[\left\lfloor \frac{\pi_*(K_X)}{n} \right \rfloor = \left
\lfloor \frac{\pi_*\pi^*(K_Y) + \pi_*(R)}{n} \right \rfloor =
K_Y + \left \lfloor \frac{\pi_*(R)}{n} \right\rfloor.\]
 Using the previous lemma we conclude that
\[\pi_*^G(\Omega_X) \cong \cO_Y\left(K_Y + \left \lfloor
\frac{\pi_*(R)}{n}\right\rfloor\right)\] and finally
\begin{eqnarray*}
\lefteqn{\textrm{dim}_k \left(H^0(X,\Omega_X)^G \right)}\\
& = & \textrm{dim}_k \left(H^0\left(Y, \pi_*^G(\Omega_X)\right)\right) \\
& = & \textrm{dim}_k
\left(H^0\left(Y, \cO_Y\left(K_Y+ \left\lfloor \frac{\pi_*(R)}{n}\right\rfloor \right) \right) \right).
\end{eqnarray*}
If $\textrm{deg}\left\lfloor \frac{\pi_*(R)}{n} \right \rfloor
= 0$ then $\left \lfloor \frac{\pi_*(R)}{n} \right\rfloor$ is
the zero divisor and we conclude that
\[\textrm{dim}_k\left(H^0(X,\Omega_X)^G\right) =
\textrm{dim}_k\left(H^0(Y, \Omega_Y)\right) = g_Y,\] as desired.
If $\textrm{deg}\left\lfloor \frac{\pi_*(R)}{n} \right \rfloor
> 0$ the divisor $K_Y + \left \lfloor \frac{\pi_*(R)}{n} \right
\rfloor$ is non-special and using the Riemann-Roch theorem (see
Theorem~1.3 on p.~295 and Example~1.3.4 on p.~296 in \cite{Ha})
we obtain
\begin{eqnarray*}
\lefteqn{\dim_k\left(H^0(X, \Omega_X)^G \right)}\\
& = & \textrm{deg}\left(K_Y + \left \lfloor \frac{\pi_*(R)}{n} \right \rfloor \right) + 1 - g_Y \\
& = & g_Y - 1 + \textrm{deg} \left \lfloor \frac{\pi_*(R)}{n} \right\rfloor ,
\end{eqnarray*}
as stated. \hfill $\Box$

{\em Proof of Theorem.} We assume that $g_X \ge 2$ and prove
conditions (i), (ii) and (iii). By replacing $G$ with the
(non-trivial) kernel $H$ of the action of $G$ on $H^0(X,
\Omega_X)$ we may assume that $G$ is non-trivial and that $G$
acts trivially on $H^0(X,\Omega_X)$: for condition~(i) we note
that if the projection $X \ra X/H$ is not tamely ramified then
the projection $\pi: X \ra Y$ cannot be tamely ramified either;
for part~(ii) we note that if the genus of $X/H$ is $0$ also
the genus of $Y=X/G$ is zero by the Hurwitz formula (see
Corollary~2.4 on p.~301 in \cite{Ha}) applied to the canonical
morphism $X/H \ra Y$; there isn't anything to
note for part~(iii) in this reduction. \\
We first suppose that $\pi$ is tamely ramified. Then we have
$R= \sum_{P \in X} (e_P-1)[P]$ by Proposition 2.2(c) on p.~300
in \cite{Ha}; hence $\left \lfloor \frac{\pi_*(R)}{n} \right
\rfloor$ is the zero divisor. Therefore we obtain
\[g_X = \textrm{dim}_k\left(H^0(X, \Omega_X)\right) =
\textrm{dim}\left(H^0(X,\Omega_X)^G\right) = g_Y\] by
Proposition~1. Substituting this equality into the Hurwitz
formula
\[2(g_X -1) = 2n (g_Y-1) + \textrm{deg}(R)\]
yields the desired contradiction because $n \ge 2$, $g_X \ge 2$
and $\textrm{deg}(R) \ge 0$. \\
We arrive at this contradiction whenever $\left \lfloor
\frac{\pi_*(R)}{n} \right \rfloor$ is the zero divisor. To
prove condition~(ii) we can therefore assume that $\left
\lfloor \frac{\pi_*(R)}{n} \right \rfloor$ is not the zero
divisor. Then Proposition~1 tells us that
\[g_X = g_Y-1 + \textrm{deg} \left \lfloor \frac{\pi_*(R)}{n}
\right \rfloor.\] Substituting this equality into the Hurwitz
formula we obtain
\[2\left(g_Y - 1 + \textrm{deg}\left \lfloor \frac{\pi_*(R)}{n}
\right \rfloor -1 \right) = 2n (g_Y -1) + \textrm{deg}(R).\]
For any $Q \in Y$ let $n_Q$ denote the coefficient of the
ramification divisor $R$ at any $P \in \pi^{-1}(Q)$ and let
$e_Q := e_P$ for any $P \in \pi^{-1}(Q)$. Rewriting the
previous equation yields
\begin{eqnarray*}
\lefteqn{(2n-2)g_Y = 2n-4 + 2 \,\textrm{deg}\left \lfloor
\frac{\pi_*(R)}{n}\right \rfloor - \textrm{deg}(R)}\\
&=& 2 \left(n-2 + \sum_{Q \in Y}
\left(\left\lfloor \frac{n}{e_Q} \frac{n_Q}{n} \right\rfloor - \frac{n}{e_Q} \frac{n_Q}{2}\right) \right)\\
&=& 2 \left(n-2 + \sum_{Q \in Y}
\left( \left\lfloor \frac{n_Q}{e_Q} \right\rfloor - \frac{n_Q}{e_Q} \frac{n}{2} \right)\right)\\
& \le & 2(n-2)
\end{eqnarray*}
because $\frac{n}{2} \ge 1$ and $\left\lfloor \frac{n_Q}{e_Q}
\right\rfloor \le \frac{n_Q}{e_Q}$ for all $Q \in Y$. Hence we
obtain $g_Y \le \frac{n-2}{n-1} < 1$ and therefore $g_Y =0$, as
stated in condition~(ii). \\
By condition~(i) the characteristic of $k$ is positive and we
may furthermore replace $G$ by a cyclic subgroup of $G$ of
order $p$ in order to prove condition~(iii). Then
condition~(iii) is the conclusion of one direction in the
following proposition. \hfill $\Box$

{\bf Proposition~2.} {\em Let $p  > 0$ and let $G$ be cyclic of
order $p$. We furthermore assume that $g_Y = 0$. Then $G$ acts
trivially on $H^0(X,\Omega_X)$ if and only if one of the
following three conditions holds:\\
(i) $p=2$.\\
(ii) $g_X =0$.\\
(iii) $p=3$ and $g_X=1$. }

{\em Proof.} Let $P_1, \ldots, P_r \in X$ be the ramified
points of $\pi: X \ra Y$ and, for $i=1, \ldots, r$, define $N_i
\in \NN$ by $\textrm{ord}_{P_i}(\sigma(\pi_i) - \pi_i) = N_i
+1$ where $\pi_i$ is a local parameter at $P_i$ and $\sigma$ is
a generator of $G$. From Lemma~1 on p.~87 in \cite{Na} we know
that $p$ does not divide $N_i$, a fact we will use several
times below. The ramification divisor $R$ of $\pi$ is equal to
$\sum_{i=1}^r(N_i+1)(p-1)[P_i]$ by Hilbert's formula for the
order of the different (see Prop.~4, \S 1, Ch.~IV on p.~72 in
\cite{Se}). Let $N:= \sum_{i=1}^r N_i$. Using the Hurwitz
formula we obtain
\[2g_X - 2 = -2p + (N+r)(p-1)\] and hence
\[\textrm{dim}_k\left(H^0(X,\Omega_X)\right) = g_X =
\frac{(N+r-2)(p-1)}{2}.\] Since $g_X \ge 0$ we obtain $r \ge
1$; that is, $\pi$ is not unramified. Therefore we have
\[\textrm{deg} \left\lfloor \frac{\pi_*(R)}{p} \right\rfloor =
\sum_{i=1}^r \left\lfloor \frac{(N_i+1)(p-1)}{p}\right\rfloor
\ge \sum_{i=1}^r \left\lfloor \frac{2(p-1)}{p}\right\rfloor = r
> 0.\] From Proposition~1 we then conclude that
\begin{eqnarray*}
\lefteqn{\textrm{dim}_k
\left(H^0(X,\Omega_X)^G\right) = g_Y -1 + \textrm{deg}
\left\lfloor \frac{\pi_*(R)}{p} \right\rfloor}\\
&=& -1 + \sum_{i=1}^r \left\lfloor \frac{(N_i+1)(p-1)}{p}\right\rfloor\\
&=& -1 + N +r + \sum_{i=1}^r \left\lfloor -\frac{N_i+1}{p}\right\rfloor.
\end{eqnarray*}
If $p=2$ the dimension of both $H^0(X,\Omega_X)$ and
$H^0(X,\Omega_X)^G$ is therefore equal to $\frac{N+r-2}{2}$. If
$g_X = 0$ both dimensions are obviously equal to~$0$. If $p=3$
and $g_X =1$ we obtain $N+r=3$ and hence $r=1$ and $N=2$; thus
both dimensions are equal to $1$. Therefore in all three of
these cases $G$ acts trivially on $H^0(X,\Omega_X)$.
This finishes the proof of one direction in Proposition~2.\\
To prove the other direction we now assume that $G$ acts
trivially on $H^0(X, \Omega_X)$ and that $p \ge 3$ and prove
that condition~(ii) or condition~(iii) holds. For each $i=1,
\ldots, r$, we write $N_i = s_i p +t_i$ with $s_i \in \NN$ and
$t_i \in \{1, \ldots, p-1\}$. We furthermore put $S:=
\sum_{i=1}^r s_i$ and $T:= \sum_{i=1}^r t_i \ge r$. Then we
have
\[ \frac{(N+r-2)(p-1)}{2} =\textrm{dim}_k(H^0(X,\Omega_X))  =
\textrm{dim}_k\left(H^0(X,\Omega_X)^G\right) = N-S-1 .\]
Rearranging this equation we obtain
\[(3-p)N - 2S = (r-2)(p-1) +2  \]
and hence
\[(-p^2 + 3p -2)S = (r-2)(p-1) +2 - (3-p)T.\]
Since $-p^2+3p-2 = - (p-1)(p-2)$ and $p \ge 3$ this equation
implies that
\[ S = \frac{(r-2)(1-p)-2 + T (3-p)}{(p-1)(p-2)}. \]
Because $S \ge 0$ the numerator of this fraction is
non-negative, that is
\begin{eqnarray*}
\lefteqn{0 \le (r-2)(1-p) - 2 + T (3-p)}\\
&\le & (r-2)(1-p) - 2 + r (3-p)\\
&=& 2 (r-1)(2-p).
\end{eqnarray*}
Hence we have $r=1$ and that numerator is $0$. We conclude that
$S=0$ and hence that $T=1$ or $p=3$. If $T=1$ we also have
$N=1$ and finally
\[g_X = \frac{(N+r-2)(p-1)}{2} = 0,\]
i.e.\ condition~(ii) holds. If $T \not=1$ and $p=3$ we obtain
$N=T=2$ and finally \[g_X = \frac{(N+r-2)(p-1)}{2} =1,\] i.e.\
condition~(iii) holds. \hfill $\Box$

{\em Acknowledgments.} The question underlying this paper goes
back to Michel Matignon. I would like to thank Niels Borne for
communicating this question to me, for sketching a proof of
condition~(i) of the above theorem and in particular for
pointing me to the above lemma.

\end{document}